\documentclass{article}

\title{\LARGE \textbf{My Research Visiting Card \\ in Hamiltonian Graph Theory}}
\author{Zh.G. Nikoghosyan\footnote{G.G. Nicoghossian (up to 1997)}\\ \\
Institute for Informatics and Automation Problems\\ National Academy of Sciences\\
P. Sevak 1, Yerevan 0014, Armenia\\ E-mail: zhora@ipia.sci.am}

\begin{document}

\maketitle

\begin{abstract}
We present eighteen exact analogs of  six well-known fundamental Theorems (due to Dirac, Nash-Williams and Jung) in hamiltonian graph theory providing alternative compositions of graph invariants. In Theorems 1-3 we give three lower bounds for the length of a longest cycle $C$ of a graph $G$ in terms of minimum degree $\delta$, connectivity $\kappa$ and parameters $\overline{p}$, $\overline{c}$ - the lengths of a longest path and longest cycle in $G\backslash C$, respectively.  These bounds have no analogs in the area involving   $\overline{p}$ and $\overline{c}$ as parameters. In Theorems 11 and 12 we give two Dirac-type results for generalized cycles including a number of fundamental results (concerning Hamilton and dominating cycles) as special cases.  Connectivity invariant $\kappa$ appears as a parameter in some fundamental results and in some their exact analogs (Theorems 3-10) in the following chronological order: 1972 (Chv\'{a}tal and Erd\"{o}s), 1981a (Nikoghosyan), 1981b (Nikoghosyan), 1985a (Nikoghosyan), 1985b (Nikoghosyan), 2000 (Nikoghosyan), 2005 (Lu, Liu, Tian), 2009 (Nikoghosyan), 2009a (Yamashita), 2009b (Yamashita), 2011a (Nikoghosyan), 2011b (Nikoghosyan).  \\


\end{abstract}

\section{The research area and the research objects}

A cycle that contains every vertex of a graph is called a Hamilton cycle, and a graph that contains a Hamilton cycle is called hamiltonian.
This terminology is used in honor of Sir William Rowan Hamilton, who, in 1865,
described an idea for a game in a letter to a friend. This game, called the icosian game,
consisted of a wooden dodecahedron  with pegs inserted at each of the
twenty vertices. These pegs supposedly represented the twenty most important cities of
the time. The object of the game was to mark a route (following the edges of the
dodecahedron) passing through each of the cities only once and finally returning to the
initial city. Hamilton sold the idea for the game. 

The classic hamiltonian problem; determining when a graph contains a Hamilton cycle, has long been fundamental in graph theory and has been intensively studied by numerous researchers. The problem of determining Hamiltonicity of a graph is one of the NP-complete problems that Karp listed in his seminal paper [6], and accordingly, one cannot hope for a simple classification of such graphs. 

Today, any path or cycle problem, related to their lengths, is really a part of this general area. However, in the last 60 years much of research in hamiltonian graph theory has focused on the following two general problems:

\begin{itemize}
\item find sufficient conditions for Hamilton and dominating cycles and for their some generalizations,
\item find lower bounds for the circumference - the length of a longest cycle in a graph.
\end{itemize}

Nowdays, we cannot even begin to completely survey the
extensive literature on hamiltonian graph theory. 

\section{Notations and definitions}

We consider only finite undirected graphs without loops or multiple edges. 
A good reference for any undefined terms is \cite{[1]}. 

A simple cycle (or just a cycle) $C$ of length $t$ is a sequence $v_1v_2...v_tv_1$ of distinct vertices $v_1,...,v_t$ with $v_iv_{i+1}\in E(G)$ for each $i\in \{1,...,t\}$, where $v_{t+1}=v_1$. When $t=2$, the cycle $C=v_1v_2v_1$ on two vertices $v_1, v_2$ coincides with the edge $v_1v_2$, and when $t=1$, the cycle $C=v_1$ coincides with the vertex $v_1$. So, by this standard definition, all vertices and edges in a graph can be considered as cycles of lengths 1 and 2, respectively.  A graph $G$ is hamiltonian if $G$ contains a Hamilton cycle, i.e. a cycle of length $n$.    If 
$Q$ is a cycle then we use  $|Q|$ to denote the length of $Q$, that is $|Q|=|V(Q)|$. We call $Q$  a large cycle in a graph $G$  if it dominates some certain subgraph structures  in $G$ in a sense that every such structure has a vertex in common with $Q$. If   $Q$  dominates all vertices in $G$ then clearly $C$ is a Hamilton cycle.  A cycle $Q$ is a dominating cycle   if it dominates all edges in $G$. A cycle $Q$ is a $PD_\lambda$ (path dominating) cycle  if it dominates all paths in $G$ of length at least some fixed integer $\lambda$. 
Finally, a cycle $Q$ is a  $QD_\lambda$ (cycle dominating) cycle  if it dominates all cycles in $G$ of length at least  $\lambda$. 

We use $n, q, \delta$ and $\alpha$ to denote the order (the number of vertices), the size (the number of edges), minimum degree and independence number of a graph. The connectivity  $\kappa$ is defined to be the minimum number of vertices whose removal disconnects $G$ or reduces it to a single vertex $K_1$.  Let $s(G)$ denote the number of components of a graph $G$. A graph $G$ is $t$-tough if $|S|\ge ts(G\backslash S)$ for every subset $S$ of the vertex set $V(G)$ with $s(G\backslash S)>1$. The toughness of $G$, denoted $\tau(G)$, is the maximum value of $t$ for which $G$ is $t$-tough (taking $\tau(K_n)=\infty$ for all $n\ge 1$). For a given longest cycle $C$ we denote by $\overline{p}$ and $\overline{c}$  the lengths of a longest path and a longest cycle in $G\backslash C$.

Let $a,b,t,k$ be integers with $k\le t$. We use $H(a,b,t,k)$ to denote the graph obtained from $tK_a+\overline{K}_t$ by taking any $k$ vertices in subgraph $\overline{K}_t$ and joining each of them to all vertices of $K_b$. Let $L_\delta$ be the graph obtained from $3K_\delta +K_1$ by taking one vertex in each of three copies of $K_\delta$ and joining them each to other. For odd $n\ge 15$, construct the graph $G_n$ from $\overline{K}_{\frac{n-1}{2}}+K_\delta+K_{\frac{n+1}{2}-\delta}$, where $n/3\le\delta\le (n-5)/2$, by joining every vertex in $K_\delta$ to all other vertices and by adding a matching between all vertices in $K_{\frac{n+1}{2}-\delta}$ and $(n+1)/2-\delta$ vertices in $\overline{K}_{\frac{n-1}{2}}$. It is easily seen that $G_n$ is 1-tough but not hamiltonian. A variation of the graph $G_n$, with $K_\delta$ replaced by $\overline{K}_\delta$ and $\delta=(n-5)/2$, will be denoted by $G^*_n$.

\section{Fundamentals and possible developments}

The earliest two theoretical results in hamiltonian graph theory  were developed in 1952 due to Dirac \cite{[3]}  in forms of a sufficient condition for a graph to be hamiltonian and a lower bound for the circumference $c$ - the length of a longest cycle in a graph, based on order $n$ and minimum degree $\delta$.\\

\noindent\textbf{Theorem A} (Dirac \cite{[3]}, 1952).  Every graph with $\delta\ge \frac{1}{2}n$ is hamiltonian.\\

\noindent\textbf{Theorem B} (Dirac \cite{[3]}, 1952).  In every 2-connected graph,  $c\ge\min\{n,2\delta\}$.\\

The exact analogs of Theorems A and B for dominating cycles was given by Nash-Williams \cite{[8]} and Jung \cite{[5]}.  \\

\noindent\textbf{Theorem C} (Nash-Williams \cite{[8]}, 1971).  Let $G$ be a 2-connected graph with $\delta\ge \frac{1}{3}(n+2)$. Then  each longest cycle in $G$ is a dominating cycle.\\

\noindent\textbf{Theorem D} (Jung \cite{[5]}, 1981).  Let $G$ be a 3-connected graph. Then either $c\ge 3\delta-3$ or each longest cycle in $G$ is a dominating cycle.\\

Theorems C and D can be easily modified to the following stronger results under additional relation $\delta\ge\alpha$ between the minimum degree $\delta$ and independence number $\alpha$.\\

\noindent\textbf{Theorem E} (Nash-Williams \cite{[8]}, 1971).  Every 2-connected graph with $\delta\ge\alpha$ and $\delta\ge\frac{1}{3}(n+2)$ is hamiltonian.\\

\noindent\textbf{Theorem F} (Jung \cite{[5]}, 1981).  In every 3-connected graph with $\delta\ge\alpha$,  $c\ge \min\left\{n,3\delta-3\right\}$.\\

In 1972, Chv\'{a}tal and Erd\"{o}s \cite{[2]} found a simple relation between connectivity $\kappa$ and independence number $\alpha$ insuring the existence of a Hamilton cycle.\\

\noindent\textbf{Theorem G} (Chv\'{a}tal and Erd\"{o}s \cite{[2]}, 1972). Every graph with $\kappa\ge\alpha$ is hamiltonian.\\

One can easily observe that the most of results  in hamiltonian graph theory are inspired by fundamental Theorems A-G
by introducing various additional new ideas, generalizations, extensions, restrictions and structural limitations such as:  

\begin{itemize}
\item \textbf{generalized and extended graph invariants} - degree sequences (P\'{o}sa type, Chv\'{a}tal type), degree sums (Ore type, Fun type), neighborhood unions, generalized degrees, local connectivity, and so on,
\item \textbf{extended list of path and cycle structures} - Hamilton, longest and dominating cycles, generalized cycles including Hamilton and dominating cycles as special cases, 2-factor, multiple Hamilton cycles, edge disjoint Hamilton cycles, powers of Hamilton cycles, $k$-ordered Hamilton cycles, arbitrary cycles, cycle systems, pancyclic-type cycle systems, cycles containing specified sets of vertices or edges,  shortest cycles, analogous path structures, and so on, 
\item \textbf{structural (descriptive) limitations} - regular, planar, bipartite, chordal and  interval graphs, graphs with forbidden subgraphs, Boolean graphs, hypercubes,  and so on,
\item\textbf{graph extensions} - hypergraphs, digraphs and orgraphs, labeled and weighted graphs, infinite graphs, random graphs, and so on.
\end{itemize}

Returning to Theorems A-F, observe that the bounds 
$$
\frac{n}{2},\ \   \frac{n+2}{3}, \ \   2\delta,   \ \   3\delta-3
$$
in these theorems are the maximum that can be press out from graph invariants $n$, $\delta$, $\alpha$ and their possible combinations. In other words, the impact of graph invariants $n$, $\delta$ and $\alpha$ on cycle structures is limited due to their nature.

A natural problem arises to find the exact analogs of fundamental Theorems A-F with alternative compositions of parameters by incorporating new advanced graph invariants  such as connectivity $\kappa$ and toughness $\tau$ - the best known graph invariants having much more impact on cycle structures in graphs.

\section{The main contribution}

In Theorems 1-3 we give three lower bounds for the length of a longest cycle $C$ of a graph $G$ in terms of minimum degree $\delta$, connectivity $\kappa$ and parameters $\overline{p}$, $\overline{c}$ - the lengths of a longest path and longest cycle in $G\backslash C$, respectively.  These bounds have no analogs in the area involving   $\overline{p}$ and $\overline{c}$ as parameters, despite the fact that the idea of using $G\backslash C$ appropriate structures lies in the base of almost all existing proof techniques in trying to construct longer cycles in graphs by the following standard procedure: choose an initial cycle
$C_0$  in $G$ and try to enlarge it by replacing a segment $P^\prime$ of $C_0$ 
with a suitable path $P^{\prime\prime}$ longer than $P^\prime$, 
having the same end vertices and passing through $G\backslash C_0$.
To find suitable $P^{\prime}$ and $P^{\prime\prime}$, one can use   
the paths or cycles (preferably large) in $G\backslash C_0$ and connections (preferably high) between
these paths (cycles) and $C_0$. The latter are closely related to $\overline{p}$,
$\overline{c}$, as well as minimum degree $\delta$ (local connections) and connectivity $\kappa$ 
(global connections). 

We give also a number of exact analogs of  fundamental Theorems A-F, concentrated especially on connectivity parameter  $\kappa$. 

The main contribution can be formulated as follows.

\begin{itemize}
\item Eighteen exact analogs of fundamental Theorems A-F were developed (Theorems 1-18) providing alternative advanced compositions of graph invariants. 
\item For fixed longest cycle $C$ in a graph $G$, the appropriate $G\backslash C$ structures  are intensively studied toward using these structures to construct longer cycles in graphs. The lengths $\overline{p}$ and $\overline{c}$ of a longest path and a longest cycle in $G\backslash C$, respectively, are incorporated into three lower bounds for $|C|$ as parameters (Theorems 1-3, having no analogs in the area), essentially improving the previous bounds.  
\item The connectivity invariant $\kappa$ (the most interesting parameter related to cycle structures in graphs) appears as a parameter in the fundamental Theorem G, as well as in a number of exact analogs (Theorems 3-10) of the fundamental Theorems A-F in the following chronological order:  1972 (Chv\'{a}tal and Erd\"{o}s \cite{[2]}), 1981a (Nikoghosyan \cite{[10]}), 1981b (Nikoghosyan \cite{[10]}), 1985a (Nikoghosyan \cite{[11]}), 1985b (Nikoghosyan \cite{[12]}), 2000 (Nikoghosyan \cite{[15]}), 2005 (Lu, Liu, Tian \cite{[7]}), 2009 (Nikoghosyan \cite{[17]}), 2009a (Yamashita \cite{[23]}), 2009b (Yamashita \cite{[23]}), 2011a (Nikoghosyan \cite{[9]}), 2011b (Nikoghosyan \cite{[9]}). 
\end{itemize}

\section{Results}

We begin with two lower bounds for the length of a longest cycle  $C$ in a graph $G$ based on parameters $\overline{p}$ and $\overline{c}$ - the lengths of a longest path and a longest cycle in $G\backslash C$, respectively. \\

\noindent\textbf{Theorem 1} \cite{[13]} (1998)  

\noindent Let $G$ be a graph and $C$ a longest cycle in $G$. Then 
$$
|C|\ge(\overline{p}+2)(\delta-\overline{p}).
$$
  
Example for sharpness: $(\kappa+1)K_{\delta-\kappa+1}+K_\kappa$.\\

\noindent\textbf{Theorem 2}  \cite{[14]} (2000)

\noindent Let $G$ be a graph and $C$ a longest cycle in $G$. Then 
$$
|C|\ge(\overline{c}+1)(\delta-\overline{c}+1).
$$

Example for sharpness: $(\kappa+1)K_{\delta-\kappa+1}+K_\kappa$.\\

For fixed $\delta$, the bound  in Theorem  2 (for example) is a strictly increasing function of  $\overline{c}$ when $0\le\overline{c}\le\delta/2$, and is a strictly decreasing function when $\delta/2\le\overline{c}\le\delta$ with global maximum
$$
\left(\left\lbrack\delta/2\right\rbrack+1\right)\left(\right\rbrack\delta/2\left\lbrack+1\right)\cong\frac{(\delta+2)^2}{4}
$$
at $\overline{c}=\lbrack\delta/2\rbrack$. 

My best result is the following lower bound for the circumference, providing increasing function on  $\overline{c}$, $\kappa$ and $\delta$.\\

\noindent\textbf{Theorem 3}  \cite{[15]} (2000)

\noindent Let $G$ be a graph with $\kappa\ge2$ and $C$ a longest cycle in $G$. If $\overline{c}\ge \kappa$ then 
$$
|C|\ge \frac{(\overline{c}+1)\kappa}{\overline{c}+\kappa+1}(\delta+2).
$$
Otherwise, 
$$
|C|\ge\frac{(\overline{c}+1)\overline{c}}{2\overline{c}+1}(\delta+2).
$$

Example for sharpness: $(\kappa+1)K_{\delta-\kappa+1}+K_\kappa$.\\

\noindent\textbf{Theorem 4}  \cite{[10]} (1981)   \  \  \     (Analog of Theorem A)

\noindent Every 2-connected graph with 
$$
\delta\ge \frac{n+\kappa}{3}
$$
is hamiltonian.\\

Examples for sharpness: $2K_\delta+K_1$; $H(1,\delta-\kappa+1,\delta,\kappa)$ $(2\le\kappa<n/2)$.\\

A short proof of Theorem 4 was given by H\"{a}ggkvist \cite{[4]}.\\

\noindent\textbf{Theorem 5}  \cite{[11]} (1985)\ \ \ (Analog of Theorem E)

\noindent Every 3-connected graph with 
$$
\delta\ge \max\left\{\frac{n+2\kappa}{4},\alpha\right\}
$$
is hamiltonian.\\

Examples for sharpness: $3K_2+K_2$; $4K_2+K_3$, $H(1,2,\kappa+1,\kappa)$.\\

The graph $4K_2+K_3$ shows that for $\kappa=3$ the minimum degree bound $(n+2\kappa)/4$ in Theorem 5 can not be replaced by $(n+2\kappa-1)/4$.\\

\noindent\textbf{Theorem 6}  \cite{[10]}  (1981)\ \ \ (Analog of Theorem B)

\noindent In every 3-connected graph,  
$$
c\ge\min\{n,3\delta-\kappa\}.
$$

Examples for sharpness: $3K_{\delta-1}+K_2$; $H(1,\delta-\kappa+1,\delta,\kappa)$.\\

\noindent\textbf{Theorem 7}  \cite{[12]} (1985)\ \ \ (Analog of Theorem F)

\noindent In every 4-connected graph with $\delta\ge\alpha$,  
$$
c\ge\min\{n,4\delta-2\kappa\}.
$$

Examples for sharpness: $4K_2+K_3$; $H(1,n-2\delta,\delta,\kappa)$; $5K_2+K_4$. The bound $4\delta-2\kappa$ in Theorem 7 is sharp for $\kappa=4$.\\

\noindent\textbf{Theorem 8} \cite{[9]} (2011) \ \ \ (Analog of Theorem F)

\noindent In every 4-connected graph with  $\delta\ge \alpha$,  
$$
c\ge\min\{n,4\delta-\kappa-4\}.
$$

Examples for sharpness: $4K_{\delta-2}+K_3$; $H(1,2,\kappa+1,\kappa)$; $H(2,n-3\delta+3,\delta-1,\kappa)$.\\

\noindent\textbf{Theorem 9}  \cite{[17]} (2011)\ \ \ (Analog of Theorem D)

\noindent If  $G$ is a 4-connected graph then either 
$$
c\ge 4\delta-2\kappa
$$
or $G$ has  a dominating cycle.    \\         

Examples for sharpness: $4K_2+K_3$; $5K_2+K_4$; $H(1,n-2\delta,\delta,\kappa)$. Theorem 9 is sharp only for $\kappa=4$ as can be seen from $5K_2+K_4$.\\

\noindent\textbf{Theorem 10}  \cite{[9]} (2011)\ \ \ (Analog of Theorem D)

\noindent If  $G$ is a 4-connected graph then either 
$$
c\ge 4\delta-\kappa-4
$$
or each longest cycle in $G$ is a dominating cycle.    \\         

Examples for sharpness: $4K_{\delta-2}+K_3$; $H(2,\delta-\kappa+1,\delta-1,\kappa)$; $H(1,2,\kappa+1,\kappa)$.\\

\noindent\textbf{Theorem 11} \cite{[16]} (2009) \ \ \ (Generalized analog of Theorems A and C)

\noindent Let $G$ be a graph, $\lambda$ a positive integer and  
$$
\delta\ge\frac{n+2}{\lambda+1}+\lambda-2.
$$
If $\kappa\ge \lambda$  then each longest cycle in $G$ is a  $CD_{\min\{\lambda,\delta-\lambda+1\}}$-cycle.\\

Examples for sharpness: $\lambda K_{\lambda+1}+K_{\lambda-1}$ $(\lambda\geq 2)$ ; $(\lambda+1)K_{\delta-\lambda+1}+K_{\lambda}$ $(\lambda\geq 1)$ ; $H(\lambda-1,\lambda,\lambda+2,\lambda+1)$ $(\lambda\geq 2)$. \\

\noindent\textbf{Theorem 12}  \cite{[16]} (2009)\ \ \ (Generalized analog of Theorems B and D)

\noindent Let $G$ be a graph and $\lambda$ a positive integer. If  $\kappa\ge \lambda+1$ then either 
$$
c\ge (\lambda+1)(\delta-\lambda+1)
$$
or each longest cycle in $G$ is a  $CD_{\min\{\lambda,\delta-\lambda\}}$-cycle. \\    

Examples for sharpness: $(\lambda+1)K_{\lambda+1}+K_{\lambda}$ $(\lambda\geq 1)$; $(\lambda+3)K_{\lambda-1}+K_{\lambda+2}$ $(\lambda\geq 2)$; $(\lambda+2)K_{\lambda}+K_{\lambda+1}$ $(\lambda\geq 1)$.\\

\noindent\textbf{Theorem 13} \cite{[16]} (2009) \ \ \ (Generalized analog of Theorem F)

\noindent Let $G$ be a graph and $\lambda$  a positive integer. If $\kappa\ge \lambda+2$  and  $\delta\ge \alpha+\lambda-1$ then  
$$
c\ge\min\{n,(\lambda+2)(\delta-\lambda)\}.
$$

Examples for sharpness: $(\lambda +2)K_{\lambda +2}+K_{\lambda +1}$; $(\lambda +4)K_{\lambda}+K_{\lambda +3}$; $(\lambda +3)K_{\lambda+1}+K_{\lambda +2}$.\\

\noindent\textbf{Theorem 14} \cite{[18]} (2011) \ \ \ (Analog of Theorem A)

\noindent Every graph with $q\le\delta^2+\delta-1$ is hamiltonian.\\

Examples for sharpness:  The bound $\delta^2+\delta-1$ in Theorem 14 can not be relaxed to $\delta^2+\delta$ since the graph $K_1+2K_\delta$ consisting of two copies of $K_{\delta+1}$ and having exactly one vertex in common, has $\delta^2+\delta$ edges and is not hamiltonian.  \\

\noindent\textbf{Theorem 15}  \cite{[19]} (2011)\ \ \ (Analog of Theorem C)

\noindent Let $G$ be a 2-connected graph. If  
$$
q\le\left\{ 
\begin{array}{lll}
8 & \mbox{when} & \mbox{ }\delta=2, \\ \frac{3(\delta-1)(\delta+2)-1}{2} & \mbox{when} & \mbox{ }%
\delta\ge3, 
\end{array}
\right. 
$$
then each longest cycle in $G$ is a dominating cycle.\\

Examples for sharpness: To show that Theorem 15 is sharp, suppose first that $\delta =2$. The graph $K_1+2K_2$ shows that the connectivity condition $\kappa\ge2$ in Theorem 15 can not be relaxed by replacing it with $\kappa\ge1$.   The graph with vertex set $\{v_1,v_2,...,v_8\}$ and edge set 
$$
\{v_1v_2,v_2v_3,v_3v_4,v_4v_5,v_5v_6,v_6v_1,v_1v_7,v_7v_8,v_8v_4\},
$$
shows that the size bound $q\le8$ can not be relaxed by replacing it with $q\le9$. Finally, the graph $K_2+3K_1$ shows that the conclusion "each longest cycle  in $G$ is a dominating cycle" can not be strengthened by replacing it with "$G$ is hamiltonian". Analogously, we can use $K_1+2K_\delta$,    $K_2+3K_{\delta-1}$  and   $K_\delta+(\delta+1)K_1$, respectively, to show that Theorem 15 is sharp when $\delta\ge3$ .  So, Theorem 15 is best possible in all respects. \\

\noindent\textbf{Theorem 16}  \cite{[22]} (2012)\ \ \ (Analog of Theorem B)

\noindent Let $G$ be a graph with $\tau>1$. Then either 
$$
c\ge \min\{n,2\delta+5\}$$
or $G$ is the Petersen graph.\\

\noindent\textbf{Theorem 17}  \cite{[20]} (2012)\ \ \ (Analog of Theorem C)

\noindent Let $G$ be a  graph with $\tau > 1$. If  
$$
\delta\ge\frac{n-2}{3}
$$
then  each longest cycle in $G$ is a dominating cycle.\\

\noindent\textbf{Theorem 18} \cite{[21]} (2012) \ \ \ (Analog of Theorem C)

\noindent Let $G$ be a  1-tough graph. If  
$$
\delta\ge\frac{n-2}{3}
$$
then  each longest cycle in $G$ is a dominating cycle unless $G$ belongs to an easily specified class of graphs.

\end{document}